\documentclass[12pt,sort&compress]{article}
\usepackage{amssymb}
\usepackage{amsfonts}
\usepackage{amsmath}
\usepackage{amssymb,amsmath}
\usepackage{amscd}
\def\a{\alpha}
\def\b{\beta}

\def\A{\mathcal {A}}

\def\b{\beta}

\def\LL{{\mathcal L}(b)}

\def\Vir{\hbox{Vir}}
\def\cl{\centerline}

\def\vs{\vspace*}
\def\W{\mathcal{W}(b)}

\def\ni{\noindent}

\def\Z{\mathbb{Z}}
\def\CC{\mathbb{C}}

\def\QED{\hfill$\Box$}

\def\pa{\partial}
\def\la{\lambda}

\numberwithin{equation}{section}
\newtheorem{theo}{Theorem}[section]
\newtheorem{defi}[theo]{Definition}
\newtheorem{coro}[theo]{Corollary}
\newtheorem{lemm}[theo]{Lemma}
\newtheorem{prop}[theo]{Proposition}
\newtheorem{clai}{Claim}
\newtheorem{rema}[theo]{Remark}

\newtheorem{case}{Case}

\makeatletter

\def\@biblabel#1{#1.~}

\makeatother

\begin{document}
\cl{\large {\bf Classification of finite irreducible conformal modules over some }}
\cl{\large {\bf Lie conformal algebras related to the Virasoro conformal algebra} \footnote{Corresponding author: Lamei Yuan(lmyuan@hit.edu.cn).}}

\vs{12pt}

\cl{Henan Wu$^{\,\dag}$,  Lamei Yuan$^{\,\ddag}$}
\cl{\small{$^{\dag}$ School of Mathematical Sciences, Shanxi University, Taiyuan 030006, China}}
\cl{\small{ $^{\ddag}$ Academy of Fundamental and Interdisciplinary
Sciences,}}\cl{\small{Harbin Institute of Technology, Harbin 150080, China}}
\cl{\small E-mails:
wuhenan@sxu.edu.cn, lmyuan@hit.edu.cn
 }\vs{6pt}
{\small\parskip .005 truein \baselineskip 3pt \lineskip 3pt

\noindent{\bf Abstract:} In this paper, we classify all finite irreducible conformal modules over a class of Lie conformal algebras $\mathcal{W}(b)$ with $b\in\CC$ related to the Virasoro conformal algebra.
 Explicitly, any finite irreducible conformal module over $\mathcal{W}(b)$ is proved to be isomorphic
to $M_{\Delta,\alpha,\beta}$ with $\Delta\neq 0$ or $\beta\neq 0$ if $b=0$, or $M_{\Delta,\alpha}$ with $\Delta\neq 0$ if $b\neq0$. As a byproduct, all finite irreducible conformal modules over the Heisenberg-Virasoro conformal algebra and the Lie conformal algebra of $\mathcal{W}(2,2)$-type are classified. Finally,
the same thing is done for the Schr\"odinger-Virasoro conformal algebra.

\vs{5pt}

\vs{5pt}

\noindent{\bf MR(2010) Subject Classification:}~ 17B10, 17B65, 17B68


\section{Introduction}

The notion of a Lie conformal algebra, introduced by Kac in \cite{K1}, encodes an axiomatic description
of the operator product expansion of chiral fields in conformal field
theory.
It has been proved to be an effective tool for the
study of infinite-dimensional Lie algebras and associative algebras (and their representations),
satisfying the sole locality property \cite{K3}.
The general structure theory, cohomology
theory and representation theory have been developed in the literatures \cite{BKV,CK,DK}. 
In particular,
all finite irreducible modules over the Virasoro conformal algebra were classified in \cite{CK} using the equivalent languages of conformal algebras and extended annihilation algebras.

In the current paper, we aim to study conformal modules over a family of Lie conformal algebras $\W$ introduced in \cite{XY} and classify all finite irreducible conformal modules over some Lie conformal algebras related to them. By definition, the Lie conformal algebra $\W$ with a parameter $b\in\CC$ is a free Lie conformal algebra generated by 
$L$ and $H$ as a $\CC[\partial]$-module, 
and satisfying 
\begin{eqnarray}\label{la-brac}
[L_\lambda L]=(\partial+2\lambda)L,\ {[L_\lambda
H]}=(\partial+(1-b)\lambda)H,\
[H_\lambda H]=0.
\end{eqnarray}
It contains the
Virasoro conformal algebra $\Vir$ as a subalgebra, which is a
free $\CC[\partial]$-module generated by $L$ such that
\begin{eqnarray}\label{Vir}
\Vir=\CC[\partial]L,\ \ \ [L_\lambda L]=(\partial+2\lambda)L.
\end{eqnarray}
Moreover, the Lie conformal algebra $\W$ has a nontrivial
abelian conformal ideal $\mathbb{C}[\partial]H$.
Thus it is neither simple nor semisimple.

 In \cite{XY}, the Lie conformal algebra $\W$ was called a $\mathcal{W}(a,b)$ Lie conformal algebra, because it
is closely related to the Lie algebra
$W(a,b)$ with $a,b\in\CC$, which is a semidirect of the centerless Virasoro algebra and a intermediate series module $A(a,b)$. In some special cases, the Lie algebras $W(a,b)$ turn out to be subalgebras of many interesting infinite-dimensional Lie algebras, such as the $\mathcal{W}$-infinity algebra $\mathcal{W}_{1+\infty}$, the Block and the
Schr\"odinger-Virasoro Lie algebras.


On the other hand, the Lie conformal algebras $\W$ with $b\in\CC$ contain some interesting special cases.
For example, $\mathcal{W}(0)$ is actually the Heisenberg-Virasoro conformal algebra, introduced in \cite{SY},
whose cohomology groups with trivial coefficients were completely determined in \cite{YW}. Whereas
$\mathcal{W}(-1)$ is isomorphic to the Lie conformal algebra of $\mathcal{W}(2,2)$ type, whose structures were studied in \cite{YW2}. Therefore, our results and methods can be directly applied to both of these two algebras. In addition,
the Heisenberg-Virasoro conformal algebra is a subalgebra of the Schr\"odinger-Virasoro conformal algebra,
which was introduced and studied in \cite{SY} leaving the problem of classifying all finite irreducible conformal modules
unsolved. In this paper, this problem is also solved. We believe that
our results and methods will be useful to the study of other non-semisimple finite conformal algebras.
This is the main motivation to present this work.


The rest of the paper is organized as follows. In Section 2, we recall the notions of Lie conformal algebra, annihilation algebra and conformal module, and related known results. At the end of this section,
we present the main results of this paper (see Theorem \ref{main}).

In Section 3, we first study the extended annihilation algebra of $\W$ and then
discuss the irreducibility property of all free nontrivial rank one $\W$-modules over $\CC[\pa]$.
Finally, we prove Theorem \ref{main} and also obtain direct and beautiful applications to the Heisenberg-Virasoro conformal algebra and the Lie conformal algebra of $\mathcal{W}(2,2)$ type. The method here is based on the techniques developed in \cite{CK}.

In Section 4, we further apply our results and methods to the Schr\"odinger-Virasoro conformal algebra and give a complete classification of all finite irreducible conformal modules over it.


Throughout this paper, all vector spaces and tensor products are over the complex field $\mathbb{C}$. In addition, we use notations $\Z$,  $\Z^+$ and $\CC^*$ to stand for the set of integers, nonnegative integers and nonzero complex numbers, respectively.

\section{Preliminaries and main results}
In this section, we recall basic definitions, notations and related results for later use, see \cite{CK, K1, K3}. Then we state our main results of this paper.
\begin{defi} \rm (\cite{K1})
A Lie conformal algebra $\mathcal {A}$ is a $\CC[\partial]$-module endowed with a $\CC$-bilinear map
$$\mathcal {A}\otimes \mathcal {A}\rightarrow \CC[\lambda]\otimes \mathcal {A},\ \  a\otimes b \mapsto [a_\lambda b],$$
called the $\la$-bracket, and
satisfying the following axioms ($a, b, c\in \mathcal {A}$),
\begin{eqnarray}
[\partial a_\lambda b]&=&-\lambda[a_\lambda b],\ \ [ a_\lambda \partial b]=(\partial+\lambda)[a_\lambda b] \ \ \mbox{(conformal\  sesquilinearity)},\label{Lc1}\\
{[a_\lambda b]} &=& -[b_{-\lambda-\partial}a] \ \ \mbox{(skew-symmetry)},\label{Lc2}\\
{[a_\lambda[b_\mu c]]}&=&[[a_\lambda b]_{\lambda+\mu
}c]+[b_\mu[a_\lambda c]]\ \ \mbox{(Jacobi \ identity)}\label{Lc3}.
\end{eqnarray}
\end{defi}

Let $\mathcal {A}$ be a Lie conformal algebra.
For each $j\in\Z^+$, we can define the {\it $j$-product} $a_{(j)}b$ of any two elements $a,b\in \A$ by the following generating series:
\begin{equation}\label{p1}
[a\,_{\lambda}\,b]=\sum_{j\in\Z^{+}}(a_{(j)}b)\frac{\lambda^j}{j!}.
\end{equation}
Then the following axioms of $j$-products hold:\vspace*{-7pt}
\begin{equation}\label{lcaj}
\aligned
&a_{(n)}b=0,\ {\rm for}\ n>>0;\\
&(\partial a)_{(n)}b=-na_{(n-1)}b;\\
&a_{(n)}b=\sum_{j\in\Z^{+}}(-1)^{n+j+1}\frac{\partial^j}{j!} b_{(n+j)}a;\\
&[a_{(m)},b_{(n)}]=\sum_{j=0}^m \begin{pmatrix}
m\\j
\end{pmatrix}(a_{(j)}b)_{(m+n-j)}.
\endaligned
\end{equation}
Actually, one can also define Lie conformal algebras using the language of $j$-products (c.f. \cite{K1}).

The annihilation algebra of a Lie conformal algebra $\mathcal {A}$ is defined as follows.
Let $\widetilde{\A}=\A\otimes\CC[t,t^{-1}]$. For any $at^m:=a\otimes t^m\in \widetilde{\A}$, define $\partial(at^m)=(\partial a)t^m$
and $\partial_t(at^m)=m at^{m-1}$, then we extend $\partial$ and $\partial_t$ linearly to the whole $\widetilde{\A}$.
Define 
\begin{equation}\label{p2}
[at^m, bt^n]=\sum_{j\in\Z^{+}}
\begin{pmatrix}
m\\j
\end{pmatrix}
(a_{(j)}b)t^{m+n-j}, \ \forall\, at^m,\, bt^n\in \widetilde{\A}.
\end{equation}
From \eqref{lcaj}, we can deduce that $\widetilde{\A}$ is a Lie algebra with respect to \eqref{p2}.
Set $\widetilde{\partial}=\partial+\partial_t$.
It can be verified that $\widetilde{\partial}\widetilde{\A}$ is an ideal of $\widetilde{\A}$.
Then
\begin{equation*}
Lie(\A)={\widetilde{\A}}/{\widetilde{\partial}\widetilde{\A}}
\end{equation*}
is also a Lie algebra.
Denote $a_{(m)}=\overline{at^m}\in Lie(\A)$.
The Lie subalgebra $\textit{Lie}(\A)^+={\rm span}\{a_{(m)}\,|\,a\in\A, m\in\Z^+\}$ is called {\it the annihilation algebra} of $\A$.
The {\it extended annihilation algebra} $\textit{Lie}(\A)^e$ is defined as the
semidirect product of the $1$-dimensional Lie subalgebra $\CC\partial$ and $\textit{Lie}(\A)^+$ with the action
$\partial(at^n)=-nat^{n-1}$.

\begin{defi}\rm (\cite{CK}) A conformal module $V$ over a Lie conformal algebra $\mathcal {A}$
is a $\mathbb{C}[\partial]$-module endowed with a $\CC$-bilinear map
$\mathcal {A}\otimes V\rightarrow
V[\lambda]$, $a\otimes v\mapsto a_\lambda v$, subject to the following relations for any $a,b\in\mathcal {A}$, $v\in V$,
\begin{eqnarray*}
&&a_\lambda(b_\mu v)-b_\mu(a_\lambda v)=[a_\lambda b]_{\lambda+\mu}v,\\
&&(\partial a)_\lambda v=-\lambda a_\lambda v,\ a_\lambda(\partial
v)=(\partial+\lambda)a_\lambda v.
\end{eqnarray*}
If $V$ is finitely generated over
$\mathbb{C}[\partial]$, then $V$ is simply called finite.
If $V$ has no nontrivial submodules, then $V$ is called irreducible.
\end{defi}

Since we only consider conformal modules, we shall simply shorten the term
``conformal module" to ``module".

Suppose that $V$ is a module over a Lie conformal algebra $\mathcal {A}$.
An element $v$ in $V$ is called an {\it invariant} if $\mathcal {A}_\la v=0$.
Let $V^0$ be the subspace of invariants of $V$.
Obviously, $V^0$ is a conformal submodule of $V$.
If $V^0=V$, then $V$ is a trivial module, i.e, a module on which $\mathcal {A}$ acts trivially.
A trivial module just admits the structure of a $\CC[\partial]$-module.
For a fixed complex constant $a$, there
is a natural trivial $\mathcal {A}$-module $\CC_a$, such that $\CC_a=\CC$ and $\pa v=a v, \mathcal {A}_\la v=0$ for $v\in\CC_a$.
The modules
$\CC_a$ with $a\in\CC$ exhaust all trivial irreducible $\mathcal {A}$-modules.
Therefore, we shall only consider \textbf{nontrivial modules} in the sequel.

The following result is due to \cite[Lemma 2.2]{K3}.
\begin{lemm}
Let $\mathcal {A}$ be a Lie conformal algebra and $V$ an $\mathcal {A}$-module.\\
1) If $\partial v= c v$ for some $c\in \CC$ and $v\in V$, then $\mathcal {A}_\la v=0$.\\
2) If $V$ is a finite module without any non-zero invariants, then $V$ is a free $\CC[\partial]$-module.
\end{lemm}

An element $v\in V$ is called a {\it torsion element} if there exists a nonzero $p(\partial)\in\CC[\partial]$ such that $p(\partial)v=0$.
For any $\CC[\partial]$-module $V$, there exists a nonzero torsion element if and only if there exists nonzero $v\in V$
such that $\partial v=c v$ for some $c\in\CC$.
A finitely generated $\CC[\partial]$-module $V$ is free if and only if $0$ is the only torsion element of $V$.

From the discussions above, we get the following result immediately.
\begin{lemm}\label{I}
Let $\mathcal {A}$ be a Lie conformal algebra.
Then every finite irreducible $\mathcal {A}$-module has no nonzero torsion element and
is free of finite rank over $\CC[\partial]$.
\end{lemm}

Assume that $V$ is a conformal module over $\A$.
We can also define {\it $j$-actions} of $\A$ on $V$ using the following generating series
\begin{equation*}
a\,_{\lambda}\,v=\sum_{j\in\Z^{+}}(a_{(j)}v)\frac{\lambda^j}{j!}.
\end{equation*}
The $j$-actions satisfy relations similar to those in (\ref{lcaj}).

Hence, if we define
\begin{equation*}
\aligned
&\overline{at^n}\cdot v=a_{(n)}v,\\
&\partial\cdot v=\partial v,
\endaligned
\end{equation*}
then $V$ can be also viewed as a module over $\textit{Lie}(\A)^e$.
Conversely, given a module $V$ over $\textit{Lie}(\A)^e$ such that for each $v\in V, a\in \A$
\begin{equation*}
\overline{at^n}\cdot v=0,\ \ n\geq N,
\end{equation*}
where $N$ is a nonnegative integer depending on $a$ and $v$.
Then $V$ can be made into a conformal $\A$-module by defining $\displaystyle a\,{}_\lambda\,v=\sum_{j\in\Z^+}\frac{\la^j}{j!}\overline{at^j}\cdot v$.

The following result presents a close connection between the representation theory of a Lie conformal algebra and that of its extended annihilation algebra.
\begin{lemm}\label{ann}
{\rm \cite{CK}}
A conformal module $V$ over a Lie conformal algebra $\A$ is precisely
a module over $\textit{Lie}(\A)^e$ satisfying the property
\begin{equation}\label{conformal}
\overline{at^n}\cdot v=0,\ \ n\geq N,
\end{equation}
for $v\in V, a\in A$, where $N$ is a nonnegative integer depending on $a$ and $v$.
\end{lemm}

\begin{rema}
By abuse of notations, we also call a Lie algebra module satisfying (\ref{conformal}) a \it conformal module over $\textit{Lie}(\A)^e$.
\end{rema}

For the Virasoro conformal algebra $\Vir$, it was shown in
\cite{CK} that all free nontrivial $\Vir$-modules of rank one
over $\mathbb{C}[\partial]$ are the following ones $(\Delta,
\alpha\in \mathbb{C})$:
\begin{eqnarray}
M_{\Delta,\alpha}=\mathbb{C}[\partial]v,\ \ L_\lambda
v=(\partial+\alpha+\Delta \lambda)v.
\end{eqnarray}
The module $M_{\Delta,\alpha}$ is irreducible if and only if
$\Delta\neq 0$. The module $M_{0,\alpha}$ contains a unique
nontrivial submodule $(\partial +\alpha)M_{0,\alpha}$ isomorphic to
$M_{1,\alpha}.$  Moreover, the modules $M_{\Delta,\alpha}$ with
$\Delta\neq 0$ exhaust all finite irreducible nontrivial
$\Vir$-modules.


The following result was given in \cite{XY}.

\begin{prop} \label{mod} All free nontrivial $\W$-modules of rank one over $\CC[\partial]$
are the following ones:
\begin{eqnarray*}
M=\mathbb{C}[\partial]v,\ L_\lambda
v=(\partial+\alpha+\Delta \lambda)v, \ H_\lambda v=\delta_{b,0}\b v, \ \mbox{for
some}\ \Delta,\a,\b\in\CC.
\end{eqnarray*}
\end{prop}

Denote the module $M$ in Proposition \ref{mod} by $M_{\Delta, \a}$ (resp. $M_{\Delta, \a,\b}$) if $b\neq0$ (resp. $b=0$).
The following theorem is our main result, which will be proved in Section 3.
\begin{theo} \label{main} Any finite irreducible $\W$-module $V$ is free of rank one over $\CC[\partial]$, and one and only one of the following
statements holds:
\begin{itemize}
\item[\rm (1)] If $b\neq 0$, $V$ is isomorphic to $M_{\Delta,\alpha}$ with $\Delta\neq 0$. 
\item[\rm (2)] If $b=0$, $V$ is isomorphic to $M_{\Delta,\alpha,\b}$ with $\Delta\neq0$ or $\b\neq 0$.
\end{itemize}
\end{theo}

\section{Finite irreducible modules over $\W$}

This section is devoted to proving the main theorem. Keeping notations in the previous section,
we first study the extended annihilation algebra $\textit{Lie}(\W)^e$ of $\W$ and then
discuss the irreducibility property
of the modules given in Proposition \ref{mod}.
Based on these works, we finally give a complete classification finite irreducible modules over $\W$.

\subsection{Extended annihilation algebra of $\W$}

In this subsection, we determine the extended annihilation algebra $\textit{Lie}(\W)^e$ of $\W$ and obtain some properties of  $\textit{Lie}(\W)^e$ for later use.

\begin{lemm}\label{A}
The annihilation algebra of $\W$ is 
 $\textit{Lie}(\W)^+=(\bigoplus_{i\geq -1}\CC L_i)\bigoplus (\bigoplus_{i\geq 0} \CC H_i)$
with
\begin{equation}\label{111}
[L_m, L_n]=(m-n)L_{m+n}, \ [L_m, H_n]=-(b m+n+b) H_{m+n}, \ [H_m, H_n]=0. 
\end{equation}
And the extended annihilation algebra $\textit{Lie}(\W)^e=\textit{Lie}(\mathcal{HV})^+ \bigoplus \CC\partial$, satisfying \eqref{111} and
\begin{equation}\label{222}
[\partial, L_n]=-(n+1)L_{n-1},\ \ [\partial, H_n]=-n H_{n-1}.
\end{equation}
\end{lemm}
\ni\ni{\it Proof.}\ \
By \eqref{la-brac} and (\ref{p1}), we have
\begin{eqnarray*}
&&L_{(0)}L=\partial L,\ L_{(1)}L=2 L, \ L_{(0)}H=\partial H,\ L_{(1)}H=(1-b)H,\\
&& L_{(j)}L= L_{(j)}H=0,\ \forall \ j\geq 2,\ \ H_{(i)}H=0,\ \forall \ i\geq 0.
\end{eqnarray*}
Then (\ref{p2}) implies  
for all $m, n\in\Z^+$ that
\begin{eqnarray*}
&&[L_{(m)}, L_{(n)}]=\sum_{j\in\Z^{+}}
\begin{pmatrix}
m\\j
\end{pmatrix}
(L_{(j)}L)_{(m+n-j)}
=(m-n)L_{(m+n-1)},\\
&&[L_{(m)}, H_{(n)}]=\sum_{j\in\Z^{+}}
\begin{pmatrix}
m\\j
\end{pmatrix}
(L_{(j)}H)_{(m+n-j)}
=-(bm+n) H_{(m+n-1)},
\end{eqnarray*}
and 
$[H_{(m)}, H_{(n)}]=0.$
Setting $L_m=L_{(m+1)}$
for $m\geq -1$ and $H_n=H_{(n)}$ for $n\in\Z^+$ gives \eqref{111}.

By definition, the extended annihilation algebra $\textit{Lie}(\W)^e$ is the
semidirect product of the $1$-dimensional Lie subalgebra $\CC\partial$ and $\textit{Lie}(\W)^+$ with
$\partial=-\partial_t$.
Thus \eqref{222} follows.
\QED\vskip6pt

The following result follows directly from \eqref{111} and \eqref{222}.
\begin{lemm}\label{center}
$\partial-L_{-1}$ is in the center of $\textit{Lie}(\W)^e$.
\end{lemm}

Denote $\LL=\textit{Lie}(\W)^e$.
Define $\LL_n=\bigoplus_{i\geq n}(\CC L_i\bigoplus\CC H_i)$.
Then we obtain a filtration of subalgebras of $\LL$:
\begin{equation*}
\LL\supset\LL_{-1}\supset\LL_0\cdots\supset\LL_n\supset\cdots,
\end{equation*}
where we use the convention that 
$H_{-1}=0$. Note that $\LL_{-1}=\textit{Lie}(\W)^+$.

The following simple facts follow directly from Lemma \ref{A}.

\begin{lemm}\label{D}
\begin{itemize}\parskip-3pt
\item[\rm(1)]
For $m, n> -1$,
$[\LL_m, \LL_n]\subset \LL_{m+n}$. And $\LL_n$ is an ideal of $\LL_0$ for all $n\in\Z^+$.
\item[\rm(2)] $[\pa, \LL_n]=\LL_{n-1}$ for all $n\in\Z^+$.
\end{itemize}
\end{lemm}


By the defining relations of $\LL$ in Lemma \ref{A}, we can obtain the following results immediately.
\begin{lemm}\label{B}
\begin{itemize}\parskip-3pt
\item[\rm(1)] If $b=0$, $[\LL_0, \LL_0]=\LL_1$.
\item[\rm(2)]  If $b\neq0$, $[\LL_0, \LL_0]=\CC H_0\bigoplus \LL_1$.
\end{itemize}
\end{lemm}

\begin{lemm}\label{E} For a fixed $N\in\Z^+$,
$\LL_0/\LL_N$ is a finite-dimensional solvable Lie algebra.
\end{lemm}
\ni\ni{\it Proof.}\ \
It suffice to show the solvability of $\LL_0/\LL_N$. Define the derived subalgebras of $\LL_0$ by
$\LL_0^{(0)}=\LL_0$ and $\LL_0^{(n+1)}=[\LL_0^{(n)}, \LL_0^{(n)}]$ for $n\in\Z^+$.

In the case of $b=0$, by Lemma {\ref{B}, $\LL_0^{(1)}=[\LL_0, \LL_0]=\LL_1$.
By induction on $n$, we have $\LL_0^{(n)}\subset\LL_n$ for all $n\in\Z^+$.

In the case of $b\neq0$, by Lemma {\ref{B}, $\LL_0^{(1)}=[\LL_0, \LL_0]=\CC H_0\bigoplus \LL_1$.
Since $[H_0, \LL_1]\subset \LL_1$, we have $\LL_0^{(2)}\subset\LL_1$.
By induction on $n$, we have $\LL_0^{(2n)}\subset\LL_n$ for all $n\in\Z^+$.

Therefore, there always exists an integer $m$ such that $\LL_0^{(m)}\subset\LL_N$ for some given $N\in\Z^+$.
This proves the lemma.
\QED

\subsection{Irreducible rank one modules over $\W$}

In this subsection, we discuss the irreducibility property of rank one $\W$-modules given in Proposition \ref{mod}.

\begin{prop}\label{mod1}
\begin{itemize}\parskip-3pt
\item[\rm(1)]
If $b\neq0$, $M_{\Delta, \a}$ is an irreducible $\W$-module if and only if $\Delta\neq 0$.
\item[\rm(2)] If $b=0$, $M_{\Delta, \a, \b}$ is an irreducible $\W$-module if and only if $\Delta\neq 0$ or $\b\neq 0$.
\end{itemize}
\end{prop}
\ni\ni{\it Proof.}\ \ (1) In this case, $H_\lambda v=0$. The irreducibility of $M_{\Delta, \a}$ as a $\W$-module is equivalent to that as a  $\Vir$-module. According to the theory of $\Vir$-modules, we get the result.

(2) In this case, $H_\lambda v=\b v$. We prove the result by discussing the possible values of $\Delta$ and $\b$.

If $\Delta\neq0$, $M_{\Delta,\alpha,\beta}$ is an irreducible $\Vir$-module.
Thus it is also an irreducible $\W$-module.

If $\Delta=0$ and $\b\neq 0$, we suppose $V$ is a nonzero submodule of $M_{0,\alpha,\beta}$.
There exists an element $u=f(\pa)v\in V$ for some nonzero $f(\pa)\in\CC[\pa]$.
If ${\rm deg}\, f(\pa)=0$, then $v\in V$ and $V=M_{0,\alpha,\beta}$.
If ${\rm deg}\, f(\pa)=m>0$, then $H_{\la}u=f(\pa+\la)(H_{\la}v)=\b f(\pa+\la)v\in V[\la]$.
The coefficient of $\la^m$ in $\b f(\pa+\la)v$ is nonzero multiple of $v$, which implies $V=M_{0,\alpha,\beta}$.
Therefore, $M_{0,\alpha,\beta}$ is irreducible.

If $\Delta=\b=0$, then $(\partial+\a)v$ generates a proper submodule of $M_{0,\alpha,0}$, which is exactly a reducible $\Vir$-module.
Thus $M_{0,\alpha,0}$ is a reducible $\W$-module. By the theory of $\Vir$-modules, $(\partial +\alpha)M_{0,\alpha,0}$ is the unique
nontrivial submodule of $M_{0,\alpha,0}$, and $(\partial +\alpha)M_{0,\alpha,0}\cong M_{1,\alpha,0}.$

This completes the proof.
\QED\vskip6pt

\subsection{Classification of finite irreducible $\W$-modules}

Since $\W$ has finite rank as a $\CC[\pa]$-module, we get a stronger property of conformal modules over
$\LL$.
\begin{lemm}\label{N}
For a conformal module $V$ over $\LL$ and an element $v\in V$, there exists an integer $N\in\Z_+$ such that
$\LL_N \cdot v=0$.
\end{lemm}
\ni\ni{\it Proof.}\ \
By Lemma \ref{ann}, $V$ is actually a module over $\LL$ satisfying
the following property:
for each $v\in V$, there exist $N_1\geq -1$ and $N_2\geq 0$ such that:
\begin{equation*}
\aligned
&L_n\cdot v=0,\ \ n\geq N_1,\\
&H_m\cdot v=0,\ \ m\geq N_2.
\endaligned
\end{equation*}
Taking $N={\rm max}\{N_1, N_2\}$ gives $L_n\cdot v=0, H_n\cdot v=0$, for all $n\geq N$.
Thus $\LL_N \cdot v=0$.
\QED\vskip6pt

Now we are in a position to prove the main Theorem of the paper. \vskip6pt

{\it Proof of Theorem \ref{main}.}\ \ \ Suppose that $V$ is an finite irreducible $\W$-module.
Let $V_n=\{v\in V\,|\,\LL_n\cdot v=0\}$.
By Lemma \ref{N}, there exists some $n$ such that $V_n$ is nontrivial, i.e., $V_n\neq \{0\}$.
Let $N$ be the minimal integer such that $V_N\neq\{0\}$ and denote $U=V_N$.
Next we proceed the proof by discussing the possible values of $N$.

If $N=-1$, then $\LL_{-1}\cdot U=0$.
Hence $U=\{v\,|\,L_\la v=H_\la v=0\}$ is a trivial $\W$-module. A contradiction.

Thus $N\geq 0$.
By \cite[Lemma 3.1]{CK}, $U$ is finite-dimensional.
By the definition of $U$, the central element $\partial-L_{-1}$ stabilize $U$.
By Schur's Lemma, 
there exists some $\a\in\CC$ such that
\begin{eqnarray}
(\partial-L_{-1})\cdot v=-\a v, \ \forall\ v\in V.
\end{eqnarray}
Thus $L_{-1}\cdot v=(\pa+\a) v$ for any $v\in V$.

If $N=0$, take a nonzero vector $v\in U$ and let $V'=\CC[\partial] v$. By Lemma \ref{I}, $V'$ is free of rank one over $\CC[\partial]$.
Obviously, $\LL_{-1}\cdot v\subset V'$.
Let $R_\partial$ (resp. $L_\partial$) be the right (resp. left) multiplication by $\partial$ in the universal enveloping algebra of $\LL$.
Using $R_\partial=L_\partial-{\rm ad}_\partial$ and the binomial formula, we have
\begin{equation*}
\LL_m \partial^k=(R_\partial)^k \LL_m=(L_\partial-{\rm ad}_\partial )^k \LL_m=\sum_{j=0}^k
 \partial^{k-j}({\rm ad}_\partial)^j\LL_m=\sum_{j=0}^{k} \partial^{k-j}\LL_{m-j}
\end{equation*}
for $m\in\Z^+$.
Hence $\LL_m \cdot V'\subset V'$ for all $m\in \Z^+$. Therefore, $V'$ is a submodule of $V$ and $V=V'$ by the simplicity of $V$.
The $\la$-action of $\W$ on $V$ is given by
\begin{equation*}
\aligned
&L_\la v=\sum_{j\in\Z^{+}}(L_{(j)}v)\frac{\lambda^j}{j!}=\sum_{j\in\Z^{+}}(L_{j-1}\cdot v)\frac{\lambda^j}{j!}=(\partial+\alpha)v,\\
&H_\la v=\sum_{j\in\Z^{+}}(H_{(j)}v)\frac{\lambda^j}{j!}=\sum_{j\in\Z^{+}}(H_{j}\cdot v)\frac{\lambda^j}{j!}=0.
\endaligned
\end{equation*}
Hence, $V\cong M_{0,\alpha}$ (resp. $V\cong M_{0,\alpha, 0}$) in the case of $b\neq 0$ (resp. $b=0$).
But  $M_{0,\alpha}$ (resp. $M_{0,\alpha, 0}$) is a reducible module over $\W$ by Proposition \ref{mod1}.
This contradicts the irreducibility of $V$.

If $N>0$, $U$ is an $\LL_0-$module since $[\LL_0, \LL_N]\subset\LL_N$ by Lemma \ref{D}.
Since $\LL_N\cdot U=0$, $U$ is actually an $\LL_0/\LL_N$-module.
Because $\LL_0/\LL_N$ is a finite-dimensional solvable Lie algebra by Lemma \ref{E},
there exists a nonzero common eigenvalue vector $u\in U$ under the action of $\LL_0/\LL_N$, and then the action of $\LL_0$.
Hence there exists a linear function $\chi$ on $\LL_0$ such that $x\cdot u=\chi(x) u$ for any $x\in\LL_0$.

Since $[\LL_0, \LL_0]$ varies because of $b$, we consider the following two cases.
\begin{case} $b\neq0$. \end{case}

In this case, $[\LL_0, \LL_0]=\CC H_0\bigoplus \LL_1$ by Lemma \ref{B}.
Then $H_0\cdot u=\chi(H_0) u=0$ and $\LL_1\cdot u=\chi(\LL_1)u=0$.
Therefore, $\chi$ is determined by $\chi(L_0)$.
We can assume $L_0\cdot u=\Delta u$ for some $\Delta\in \CC$.
Then $\Delta\neq 0$ (otherwise, $N=0$).
Let $V'=\CC[\partial]u$.
Then $\LL_1\cdot u\subset V'$.
By similar discussions as in the case of $N=0$, $V'$ is a submodule of $V$. Thus $V'=V$.
The $\la$-action of $\W$ on $V$ is given by $H_\la u=0$ and
\begin{equation*}
L_\la u=\sum_{j\in\Z^{+}}(L_{(j)}u)\frac{\lambda^j}{j!}=\sum_{j\in\Z^{+}}(L_{j-1}\cdot u)\frac{\lambda^j}{j!}=L_{-1}\cdot u+(L_0\cdot u)\la=(\partial+\alpha+\Delta\lambda)u.
\end{equation*}
Hence, $V\cong M_{\Delta,\alpha}$, which is irreducible by Proposition \ref{mod1} (1).

\begin{case} $b=0$.
\end{case}

In this case, $[\LL_0, \LL_0]=\LL_1$ by Lemma \ref{B}.
Then $L_0\cdot u=\chi(L_0)u$, $H_0\cdot u=\chi(H_0) u$ and $\LL_1\cdot u=\chi(\LL_1)u=0$.
We can assume $\chi(L_0)=\Delta$ and $\chi(H_0)=\b$ for some $\Delta, \b\in \CC$.
Let $V'=\CC[\partial]u$.
Then $\LL_{-1}\cdot u\subset V'$ and $V'$ is a submodule of $V$. Thus $V'=V$.
The $\la$-action of $\W$ on $V$ is given by 
\begin{equation*}
\aligned
&L_\la u=\sum_{j\in\Z^{+}}(L_{(j)}u)\frac{\lambda^j}{j!}=\sum_{j\in\Z^{+}}(L_{j-1}\cdot u)\frac{\lambda^j}{j!}=L_{-1}\cdot u+(L_0\cdot u)\la=(\partial+\alpha+\Delta\lambda)u.\\
&H_\la u=\sum_{j\in\Z^{+}}(H_{(j)}u)\frac{\lambda^j}{j!}=\sum_{j\in\Z^{+}}(H_{j}\cdot u)\frac{\lambda^j}{j!}=H_0\cdot u=\b u.
\endaligned
\end{equation*}
Hence, $V\cong M_{\Delta,\alpha,\beta}$.
By Proposition \ref{mod1} (2), $M_{\Delta,\alpha,\beta}$ is irreducible if and only if $\Delta\neq 0$ or $\b\neq 0$.

This completes the proof.\QED\vskip6pt

By Theorem \ref{main}, the following results are immediate.

\begin{coro} Any finite irreducible
 module of the Heisenberg-Virasoro conformal algebra over $\CC[\partial]$ is
of the form ($\Delta, \alpha,\beta\in\CC$)
\begin{eqnarray*}
M_{\Delta,\alpha,\beta}=\mathbb{C}[\partial]v,\ L_\lambda
v=(\partial+\alpha+\Delta \lambda)v, \ H_\lambda v=\b v,
\end{eqnarray*}
where $\Delta\neq 0$ or $\beta\neq 0$.
\end{coro}
\begin{coro}
Any finite irreducible module of the Lie conformal algebra of $\mathcal{W}(2,2)$-type over $\CC[\partial]$ is
of the form 
\begin{eqnarray*}
M_{\Delta,\alpha}=\mathbb{C}[\partial]v,\ L_\lambda
v=(\partial+\alpha+\Delta \lambda)v, \ H_\lambda v=0, \ \Delta\in\CC^*, \ \alpha\in\CC.
\end{eqnarray*}¡¡¡¡
\end{coro}

\section{Classification of finite irreducible modules over the Schr\"odinger-Virasoro conformal algebra}
\setcounter{case}{0}

In this section, we apply the results and methods in Section 3 to the Schr\"odinger-Virasoro conformal algebra, and classify all the finite irreducible conformal modules over it.

The Schr\"odinger-Virasoro conformal algebra, introduced in \cite{SY},
is 
a finite Lie conformal algebra
$\mathcal{SV}=\mathbb{C}[\partial]L\bigoplus
\mathbb{C}[\partial]M\bigoplus\mathbb{C}[\partial]Y$, endowed with the
following  nontrivial $\lambda$-brackets\vspace*{-7pt}
\begin{eqnarray}
&&[L_\lambda L]=(\partial+2\lambda)L,\ \ \ [Y_\la
L]=(\frac12\pa+\frac32\la)Y, \label{lamda-bracket11}\\
&&{[L_\lambda Y]}=(\partial+\frac32\lambda)Y,\ \ \,{[Y_\lambda
Y]}=(\partial+2\lambda)M,\label{lamda-bracket22}\\[3pt]
&&{[L_\lambda M]}=(\partial+\lambda)M,\ \ \,[M_\la L]=\la
M.\label{lamda-bracket33}
\end{eqnarray}
Note that $\mathcal{SV}$ contains the Heisenberg-Virasoro conformal algebra $\mathbb{C}[\partial]L\bigoplus
\mathbb{C}[\partial]M\cong\mathcal{W}(0)$ as a subalgebra.
Moreover, $\mathbb{C}[\partial]M$ is an abelian ideal of $\mathcal{SV}$ and the corresponding quotient $\mathcal{SV}/\mathbb{C}[\partial]M$ is isomorphic to $\mathcal{W}(-1/2)$.

The maximal
formal distribution Lie algebra corresponding to the
Schr\"odinger-Virasoro conformal algebra $\mathcal{SV}$ is the
Schr\"odinger-Virasoro Lie algebra, which was introduced
in the context of non-equilibrium statistical physics during the
process of investigating the free Schr\"{o}dinger equations
\cite{H1,H3}. By definition, the
Schr\"odinger-Virasoro Lie algebra is an
infinite-dimensional Lie algebra with a $\CC$-basis $\{L_n,
M_n,Y_p\,, \big|\,n\in \Z,p\in \frac{1}{2}+\Z\}$ and satisfying the following
non-vanishing Lie bracket\vspace*{-5pt}s ($m,n\in\Z$ and $p,q\in\frac12+\Z$)
\begin{eqnarray}\label{LB}\begin{array}{lll}
&&[L_m,L_{n}]=(m-n)L_{m+n},\ \ \
[L_m,M_n]=-nM_{m+n},\\[6pt]
&&[\,Y_p\,,Y_{q}\,]=(p-q)M_{p+q},\ \ \ \ \
[\,L_m,Y_p\,]=(\frac{m}{2}-p)Y_{m+p}.\end{array}
\end{eqnarray}

The corresponding annihilation algebra in this case is proved to be 
\begin{eqnarray}
\textit{Lie}(\mathcal{SV})^+= \sum_{m\geq -1}\CC L_m+\sum_{n\geq 0}\CC M_n +\sum_{p\geq -\frac{1}{2}}\CC Y_{p},
\end{eqnarray}
which is a subalgebra of the Schr\"odinger-Virasoro Lie algebra.
And the the extended annihilation algebra is $\textit{Lie}(\mathcal{SV})^e=\textit{Lie}(\mathcal{SV})^+ \bigoplus \CC\partial$ with
\begin{equation}\label{LB2}
[\partial, L_m]=-(m+1)L_{m-1}, \
[\partial, M_n]=-nM_{n-1},\ [\partial, Y_p]=-(p+\frac12)Y_{p-1}.
\end{equation}
As in the previous paragraph, we have a filtration for $\textit{Lie}(\mathcal{SV})^e$ of the form
\begin{equation*}
\mathcal{L}_{-1}\supset\mathcal{L}_0\cdots\supset\mathcal{L}_n\supset\cdots
\end{equation*}
where $\mathcal{L}_n=\bigoplus\limits_{i\geq n}(\CC L_i\oplus\CC M_i\oplus\CC Y_{i+\frac12})$
with the convention $M_{-1}=0$.
Similar to the arguments in Subsection 3.1, we have $[\mathcal{L}_0, \mathcal{L}_0]=\CC Y_{\frac{1}{2}}+\mathcal{L}_1$ and $[\partial, \mathcal{L}_k]=\mathcal{L}_{k-1}$ for all $k\geq0$.
Moreover, $\pa-L_{-1}$ is also a central element by (\ref{LB}) and (\ref{LB2}).

The following result is due to \cite{SY}.

\begin{prop}\label{p3}
All free nontrivial $\mathcal{SV}$-modules of rank one over $\mathbb{C}[\pa]$
are as follows,
\begin{eqnarray*}
M_{\Delta,\alpha}=\mathbb{C}[\partial]v,\ L_\lambda
v=(\partial+\alpha+\Delta \lambda)v, \ M_\lambda v=Y_\lambda v=0,\
\mbox{for some}\ \Delta,\a\in\CC.
\end{eqnarray*}
The module $M_{\Delta,\alpha}$ is irreducible if and only if $\Delta\neq 0$.
\end{prop}

 The following result shows that all finite irreducible $\mathcal{SV}$-modules are free of rank one and thus of the kind in Proposition \ref{p3}.
\begin{theo} All finite irreducible $\mathcal{SV}$-modules are of the following kind
\begin{eqnarray}\label{sv}
M_{\Delta,\alpha}=\mathbb{C}[\partial]v,\ L_\lambda
v=(\partial+\alpha+\Delta \lambda)v, \ M_\lambda v=Y_\lambda v=0,\
\mbox{for some}\ \Delta\in\CC^*,\a\in\CC.
\end{eqnarray}
\end{theo}
\ni\ni{\it Proof.}\ \
Suppose that $V$ is a nontrivial finite irreducible $\mathcal{SV}$-module.
It is also a conformal module over $\textit{Lie}(\mathcal{SV})^e$ by Lemma \ref{ann}.
Using similar arguments as in the proof of Theorem \ref{main}, we can find a nonzero vector
$v$ such that
\begin{equation*}
L_i\cdot v=M_i\cdot v=Y_{i+\frac{1}{2}}\cdot v=0, \ \ \ \forall\,i\in\Z^+\setminus\{0\},
\end{equation*}
and $L_{0}\cdot v=\Delta v, M_{0}\cdot v=\b v$ for some $\Delta, \b\in\CC$. 
Since $L_{-1}-\pa$ is a central element of $\textit{Lie}(\mathcal{SV})^e$, we can also prove that $L_{-1}-\pa$ acts on $V$ by some scalar $\a\in\CC$.

We continue the proof by discussing the possible values of $\b$.

\begin{case}
$\b=0$
\end{case}

In this case, $M_{\la}v=0$.
Since $\CC[\pa]M$ is an ideal of $\mathcal{SV}$ and $V$ is irreducible, $M_{\la}V=0$.
Then the irreducibility of $V$ as a  $\mathcal{SV}$-module is equivalent to that of $V$ as a $\mathcal{SV}/\CC[\pa]M\cong \mathcal{W}(-1/2)$-module.
By Theorem \ref{main}, $V$ is freely generated by a vector $u$ and
\begin{eqnarray*}
L_\lambda u=(\partial+\alpha+\Delta \lambda)u,\ Y_\lambda u=0.
\end{eqnarray*}
for some $\Delta\in\CC^*$.
Then $V$ is of the kind in (\ref{sv}).




\begin{case}
$\b\neq0$
\end{case}

Set
\begin{equation*}
\mathcal{Y}=\textrm{span}_\CC\{Y_{-\frac{1}{2}},L_{-1},\pa\}.
\end{equation*}
Then $\textit{Lie}(\mathcal{SV})^e$ has a decomposition of vector spaces $\textit{Lie}(\mathcal{SV})^e=\mathcal{Y}\bigoplus \mathcal{L}_0$. By PBW Theorem, the universal envelopping algebra of $\textit{Lie}(\mathcal{SV})^e$ is $U(\textit{Lie}(\mathcal{SV})^e)=U(\mathcal{Y})\bigotimes U(\mathcal{L}_0)$, where $U(\mathcal{Y})=\textrm{span}_\CC\{L_{-1}^i\pa^j Y_{-\frac{1}{2}}^k\,|\,i, j, k\in\Z^+\}$, as a vector space over $\CC$.
Since $V$ is an irreducible module over $\textit{Lie}(\mathcal{SV})^e$, $V$ is generated by $v$.
Then we have
\begin{equation*}
V=U(\textit{Lie}(\mathcal{SV})^e)\cdot v=U(\mathcal{Y})\cdot v=\sum_{i,j\in\Z^+}\CC\pa^i Y_{-\frac{1}{2}}^j\cdot v.
\end{equation*}
Let $v_j=Y_{-\frac{1}{2}}^j\cdot v$ for any $j\in\Z^+$, then $V=\sum_{j\in\Z^+}\CC[\pa] v_j$.

\begin{clai}
There exists a nonzero polynomial $m(x)\in\CC[x]$ such that
\begin{equation}\label{sv1}
m(Y_{-\frac{1}{2}})\cdot v=0.
\end{equation}
\end{clai}

Since $V$ is a free $\CC[\pa]$-module of finite rank,
there exists $m\in\Z^+$ and $f_1(\pa),\ldots,f_m(\pa)\in\CC[\pa]$ such that
\begin{equation}\label{sv2}
\sum_{j=0}^m f_j(\pa)\cdot v_j=0.
\end{equation}
Denote $n={\rm max}\{{\rm deg\,} f_1(\pa),\ldots, {\rm deg\,}f_m(\pa)\}$.
Using $R_\partial=L_\partial-{\rm ad}_\partial$ and the binomial formula, we have for $k\leq n$
\begin{equation*}
M_n\pa^k=(L_\partial-{\rm ad}_\partial)^kM_n=\sum_{i=0}^k a_i\pa^{k-i}M_{n-i},
\end{equation*}
where
\begin{equation*}
a_i=\begin{pmatrix}
k\\i
\end{pmatrix}
\begin{pmatrix}
n\\i
\end{pmatrix}
i!\in\CC^*.
\end{equation*}
By action of $M_n$ on (\ref{sv2}) and noting that
\begin{equation*}
M_i\cdot v=
\begin{cases}
0,\ \  i>0\\
\b v,\ \  i=0
\end{cases}
\end{equation*}
we get the desired result.

Now let $m(x)$ be a nonzero polynomial with minimal degree ($\geq 1$) such that $m(Y_{-\frac{1}{2}})\cdot v=0$.
Let $Y_{\frac{1}{2}}$ act on (\ref{sv1}), we have
\begin{equation*}
Y_{\frac{1}{2}}m(Y_{-\frac{1}{2}})\cdot v=[Y_{\frac{1}{2}}, m(Y_{-\frac{1}{2}})]\cdot v=m'(Y_{-\frac{1}{2}})[Y_{\frac{1}{2}}, Y_{-\frac{1}{2}}]\cdot v=\b m'(Y_{-\frac{1}{2}})\cdot v=0,
\end{equation*}
which gives $m'(Y_{-\frac{1}{2}})\cdot v=0$ since $\b\neq 0$.
A contradiction.

Thus $\b=0$ and the result is proved.
\QED

\vspace{4mm} \noindent\bf{\footnotesize Acknowledgements.}\ \rm
{\footnotesize This work was supported by National Natural Science
Foundation grants of China (11301109, 11526125) and the Research Fund for the Doctoral Program of Higher Education (20132302120042).}\\
\vskip18pt \small\footnotesize
\parskip0pt\lineskip1pt


\parskip=0pt\baselineskip=1pt
\begin{thebibliography}{9999}
\def\RE#1{\bibitem{#1}\label{#1}}

\bibitem{BKV} Bakalov B., Kac V., Voronov A., Cohomology of
conformal algebras, {\it Comm. Math. Phys.}, {\bf 200} (1999)
561--598.

\bibitem {CK} Cheng S.-J., Kac V., Conformal modules, {\it Asian
J. Math.}, {\bf 1}(1) (1997) 181--193.

\bibitem {DK} D'Andrea A., Kac V., Structure theory of finite
conformal algebras, {\it Sel. Math., New Ser.}, {\bf 4} (1998)
377--418.

\bibitem{H1} Henkel M., Schr\"{o}dinger invariance and strongly
anisotropic critical systems, {\it J. Stat. Phys.}, {\bf 75} (1994) 1023--1029.

\bibitem{H3}Henkel M., Unterberger J., Schr\"{o}dinger invariance and
space-time symmetries, {\it Nuclear Physics B}, {\bf 660} (2003) 407--412.

\bibitem {K1} Kac V.,  Vertex algebras for beginners. Univ. Lect. Series 10, AMS
(1996). Second edition 1998.


\bibitem {K3} Kac V., The idea of locality, in: H.-D. Doebner, et al. (Eds.),
Physical Applications and Mathematical Aspects of Geometry, Groups
and Algebras, World Sci. Publ., Singapore, (1997) 16--32.


\bibitem{SY} Su Y., Yuan L., Schr\"{o}dingger-Virasoro Lie conformal
algebra, {\it J. Math. Phys.}, {\bf 54} (2013) 053503, 16pp.

\bibitem{XY} Xu Y., Yue X., $W(a,b)$ Lie conformal algebra and its conformal
module of rank one, {\it Algebra Colloq.}, {\bf 22}(3) (2015) 405--412.

\bibitem{YW} Yuan L., Wu H., Cohomology of Heisenberg-Virasoro conformal algebra, to appear in {\it J. Lie Theory}.

\bibitem{YW2} Yuan L., Wu H., Structures of a Lie conformal algebra of $W(2,2)$-type, arXiv:1601.06588.
\end{thebibliography}
\end{document}